\documentclass{article}

\usepackage[cmtip,arrow]{xy}
\usepackage{amsmath,amssymb,pb-diagram,pb-xy}
\usepackage{palatino, url, multicol}
\def \a{{\mathfrak a}}
\def \al{\alpha}

\def \Aut{{\rm Aut}}
\def \BGL{{\rm BGL}}

\def \CC{{\cal C}}
\def \CE{{\cal E}}
\def \CF{{\cal F}}

\def \CO{{\cal O}}
\def \CP{{\cal P}}

\def \coker{{\rm coker}}
\def \coim{{\rm coim}}
\def \df{\ \begin{array}{c} _{\rm def}\\ ^{\displaystyle =}\end{array}\ }

\def \F{{\mathbb F}}

\def \GL{{\rm GL}}
\def \Hom{{\rm Hom}}
\def \id{{\rm id}}
\def \im{{\rm im}}

\def \Mor{{\rm Mor}}

\def \N{{\mathbb N}}
\def \onto{\to\hspace{-8pt}\to}
\def \ph{\varphi}

\def \prf{{\bf Proof: }}
\def \p{{\mathfrak p}}
\def \Per{{\rm Per}}

\def \QB{{Q\CC}}
\def \Quot{{\rm Quot}}
\def \qed{\ifmmode\eqno \square 
		\else\noproof\vskip 12pt plus 3pt minus 9pt \fi}
\def \noproof{{\unskip\nobreak\hfill\penalty50\hskip2em		\hbox{}%
     \nobreak\hfill $\square$\parfillskip=0pt%
     \finalhyphendemerits=0\par}}

\def \Sets{{\rm Sets}}
\def \setminus{\begin{picture}(18,10)\put(4,6)
                {\line(2,-1){10}}\end{picture}}

\def \spec{{\rm spec}}

\def \Z{{\mathbb Z}}
\def \({\left(}
\def \){\right)}
\def \={{\ =\ }}

\renewcommand{\sp}[1]{\left\langle #1\right\rangle}

\newtheorem{theorem}{Theorem}

\newtheorem{lemma}[theorem]{Lemma}

\begin{document}

\pagestyle{myheadings} \markright{ZETA FUNCTIONS OVER $\F_1$}

\title{Remarks on  zeta functions and K-theory\\ over\\ 
$\F_1$}
\author{Anton Deitmar}
\date{}
\maketitle

{\bf Abstract:} We show that the notion of zeta functions over the field of one element $\F_1$, as given in special cases by Soul\'e, extends naturally to all $\F_1$-schemes as defined by the author in an earlier paper. We further give two constructions of K-theory for affine schemes or $\F_1$-rings, we show that these coincide in the group case, but not in general.

$ $

\begin{multicols}{2}
\tableofcontents
\section*{Introduction}
Soul\'e \cite{Sou}, inspired by Manin \cite{Ma}, gave a definition of zeta functions over the field of one element $\F_1$. We describe this definition as follows.
Let $X$ be a scheme of finite type over $\Z$.
For a prime number $p$ one sets after Weil,
$$
Z_X(p,T)\df \exp\( \sum_{n=1}^\infty \frac{T^n}n \# X(\F_{p^n})\),
$$
where $\F_{p^n}$ denotes the field of $p^n$ elements.
This is the local zeta function over $p$, and the global zeta function of $X$ is given as
$$
\zeta_{X|\Z}(s)\df \prod_p Z_X(p,p^{-s})^{-1}.
$$

Soul\'e considered in \cite{Sou} the following condition: Suppose there exists a polynomial $N(x)$ with integer coefficients such that $\# X(\F_{p^n})=N(p^n)$ for every prime $p$ and every $n\in\N$.
Then $Z_X(p,p^{-s})^{-1}$ is a rational function in $p$ and $p^{-s}$.
The vanishing order at $p=1$ is $N(1)$.
One may thus define
$$
\zeta_{X|\F_1}(s)\= \lim_{p\to 1}\frac{Z_X(p,p^{-s})^{-1}}{(p-1)^{N(1)}}.
$$
One computes that if $N(x)=a_0+a_1x+\dots +a_nx^n$, then
$$
\zeta_{X|\F_1}(s)\= s^{a_0}(s-1)^{a_1}\cdots (s-n)^{a_n}.
$$

In the paper \cite{F1} there is given a definition of a scheme over $\F_1$ as well as an ascent functor $\cdot\otimes\Z$ from $\F_1$-schemes to $\Z$-schemes.
We say that a $\Z$-scheme is \emph{defined over $\F_1$}, if it comes by ascent from a scheme over $\F_1$.
The natural question arising is whether schemes defined over $\F_1$ satisfy Soul\'e's condition.

Simple examples show that this is not the case.
However, schemes defined over $\F_1$ satisfy a slightly weaker condition which serves the purpose of defining $\F_1$-zeta functions as well, and which we give in the following theorem.

\begin{theorem}\label{zetapoly}
Let $X$ be a $\Z$-scheme defined over $\F_1$.
Then there exists a natural number $e$ and a polynomial $N(x)$ with integer coefficients such that
for every prime power $q$ one has
$$
(q-1,e)=1\ \ \Rightarrow\ \ \# X_\Z(\F_{q})\= N(q).
$$
This condition determines the polynomial $N$ uniquely (independent of the choice of $e$).
We call it the \emph{zeta-polynomial} of $X$.
\end{theorem}

With this theorem, we can define the zeta function of an arbitrary $\F_1$-scheme $X$ as
$$
\zeta_{X|\F_1}(s)\= s^{a_0}(s-1)^{a_1}\cdots (s-n)^{a_n},
$$
if $N_X(x)=a_0+a_1x+\cdots a_n x^n$ is its zeta-polynomial.

We also define its \emph{Euler characteristic} as
$$
\chi(X)\= N_X(1)\= a_1+\cdots +a_n.
$$
This definition is due to Soul\'e \cite{Sou}. 
We repeat the justification, which is based on the Weil conjectures.

Suppose that $X/\F_p=X_Z\times_\Z\F_p$ is a smooth projective variety over the finite field $\F_p$.
Then the Weil conjectures, as proven by Deligne, say that
$$
Z_{X_\Z}(p,T)\=\prod_{l=0}^m P_l(T)^{(-1)^{l+1}},
$$
with
$$
P_l(T)\=\prod_{j=1}^{b_l}(1-\al_{l,j}T),
$$
satisfying $|\al_{l,j}|=p^{l/2}$, where $b_l$ is the $l$-th Betti-number.

On the other hand, suppose that $\# X(\F_{p^n})=N(p^n)$ holds for every $n\in\N$, where $N(x)=a_0+a_1x+\dots +a_nx^n$ is the zeta-polynomial, then one gets
$$
Z_{X_\Z}(p,T)\= \prod_{k=0}^n(1-p^{k}T)^{-a_k}.
$$
Comparing these two expressions, one gets
$$
b_l\=\begin{cases} a_{l/2} & l\ {\rm even},\\
0 & l\ {\rm odd}.\end{cases}
$$
So $\sum_{k=0}^na_k\=\sum_{l=0}^m (-1)^lb_l$ is the Euler characteristic.

For explicit computations of zeta functions and Euler numbers over $\F_1$, see \cite{Kuro}.

Next for K-theory.
Based on the idea of Tits, that $\GL_n(\F_1)$ should be the permutation group $\Per(n)$, Soul\'e also suggested that
$$
K_i(\F_1)\= \pi_i(B(\Per(\infty))^+),
$$
which is known to coincide with the stable homotopy group of the spheres, $\pi_i^s\=\lim_{k\to\infty}\pi_{i+k}(S^k)$. (The $+$ refers to Quillen's $+$ construction.)
More general, for a monoid $A$, or an $\F_1$-ring $\F_A$, one has
$$
\GL_n(A)\=\GL_n(\F_A)\= A^n\rtimes\Per(n).
$$
Setting $\GL(A)=\lim_{n\to\infty}\GL_n(A)$, one lets
$$
K_i^+(A)\= \pi_i(BGL(A)^+).
$$
On the other hand, one considers the category $\CP$ of all finitely generated projective modules over $A$ and defines
$$
K_i^Q(A)\= \pi_{i+1}(BQ\CP),
$$
where $Q$ means Quillen's $Q$-construction.
It turns out that $\pi_1(BQ\CP)$ coincides with the Grothendieck group $K_0(\CP)$ of $\CP$.
If $A$ is a group, these two definitions of K-theory agree, but not in general.

A calculation shows, that if $A$ is an abelian group, then
$$
K_i(A)\= \begin{cases} \Z\times A & i=0,\\ \pi_i^s & i>0.
\end{cases}
$$
So, for general $A$, since one has $K^+(A)=K^+(A^\times)$, this identity completely computes $K^+$.
Furthermore, for every $A$ one has a canonical homomorphism $K_i^+(A)\to K_i^Q(A)$.

I thank Jeff Lagarias for his remarks on an earlier version of this paper.

\section{$\F_1$-schemes}
For basics on $\F_1$-schemes we refer to \cite{F1}.

In this paper, a ring will always be commutative with unit and a monoid will always be commutative.
An \emph{ideal} $\a$ of a monoid $A$ is a subset with $A\a\subset \a$.
A \emph{prime ideal} is an ideal $\p$ such that $S_\p=A\setminus\p$ is a submonoid of $A$.
For a prime ideal $\p$ let $A_\p=S_\p^{-1}A$ be the \emph{localisation} at $\p$.
The \emph{spectrum} of a monoid $A$ is the set of all prime ideals with the obvious Zariski-topology (see \cite{F1}).
Similar to the theory of rings, one defines a structure sheaf $\CO_X$ on $X=\spec(A)$, and one defines a \emph{scheme over $\F_1$} to be a topological space together with a sheaf of monoids, locally isomorphic to spectra of monoids.

A $\F_1$-scheme $X$ is \emph{of finite type}, if it has a finite covering by affine schemes $U_i=\spec(A_i)$ such that each $A_i$ is finitely generated.

For a monoid $A$ we let $A\otimes\Z$ be the monoidal ring $\Z[A]$.
This defines a functor from monoids to rings which is left adjoint to the forgetful functor that sends a ring $R$ to the multiplicative monoid $(R,\times)$.
This construction is compatible with gluing, so one gets a functor $X\mapsto X_\Z$ from $\F_1$-schemes to $\Z$-schemes.

\begin{lemma}
$X$ is of finite type if and only if $X_\Z$ is a $\Z$-scheme of finite type.
\end{lemma}

\prf
If $X$ is of finite type, it is covered by finitely many affines $\spec(A_i)$, where $A_i$ is finitely generated, hence $\Z[A_i]$ is finitely generated as a $\Z$-algebra and so it follows that $X_\Z$ is of finite type.

Now suppose that $X_\Z$ is of finite type.
Consider a covering of $X$ by open sets of the form $U_i=\spec(A_i)$.
then one gets an open covering of $X_\Z$ by sets of the form $\spec(\Z[A_i])$, with the spectrum in the ring-sense.
Since $X_\Z$ is compact, we may assume this covering finite.
As $X_\Z$ is of finite type, each $\Z[A_i]$ is a finitely generated $\Z$-algebra.
Let $S$ be a generating set of $A_i$. Then it generates $\Z[A_i]$, and hence it contains a finite generating set $T$ of $\Z[A_i]$.
Then $T$ also generates $A_i$ as a monoid, so $A_i$ is finitely generated.
\qed

\section{Proof of Theorem \ref{zetapoly}}\label{pfthm1}
We will show uniqueness first.

\begin{lemma}
For every natural number $e$ there are infinitely many prime powers $q$ with $(q-1,e)=1$.
\end{lemma}

\prf
Write $e=2^km$ where $m$ is odd.
Let $n\in\N$.
The number $2^n$ is a unit modulo $m$ and hence there are infinitely many $n$ such that $2^n\equiv 1 $ modulo $m$.
Replacing $n$ by $n+1$ we see that there are infinitely many $n$ such that $2^n\equiv 2$ modulo $m$ and hence $2^n-1\equiv 1$ modulo $m$.
As $2^n-1$ is odd, it follows $(2^n-1,e)=1$ for every such $n$.
\qed

Now for the uniqueness of $N$.
Suppose that the pairs $(e,N)$ and $(e',N')$ both satisfy the theorem. Then for every prime power $q$ one has
$$
(q-1,ee')\= 1\ \ \Rightarrow\ \ N(q)\=\# X(\F_q)\=N'(q).
$$
As there are infinitely many such prime powers $q$, it follows that $N(x)=N'(x)$, as claimed.

We start on the existence of $N$.
For a finite abelian group $E$ define its \emph{exponent} $m=\exp(E)$ to be the smallest number $m$ such that $x^m=1$ for every $x\in G$.
The exponent is the least common multiple of the orders of elements of $G$.
A finitely generated abelian group $G$ is of the form $\Z^r\times E$ for a finite group $E$.
Then $r$ is called the \emph{rank} of $G$ and the exponent of $E$ is called the \emph{exponent} of $G$.

For a finitely generated monoid $A$ we denote by $\Quot(A)$ its quotient group.
This group comes about by inverting every element in $A$.
It has a natural morphism $A\to\Quot(A)$ and the universal property that every morphism from $A$ to a group factorizes uniquely over $A\to\Quot(A)$.
In the language of \cite{F1}, $\Quot(A)$ coincides with the stalk $\CO_\eta=A_\eta$ at the generic point $\eta$ of $\spec(A)$.

We define the \emph{rank} and \emph{exponent} to be the rank and exponent of $\Quot(A)$.
Note that for a finitely generated monoid $A$ the spectrum $\spec(A)$ is a finite set.
Hence the underlying space of a scheme $X$ over $\F_1$ of finite type is a finite set.
We then define the exponent of $X$ to be the least common multiple of the numbers $\exp(\CO_\p)$, where $\p$ runs through the finite set $X$.

Let $X$ be a scheme over $\F_1$ of finite type.
We may assume that $X$ is connected.
Let $e$ be its exponent.
Let $q$ be a prime power and let $D_q$ be the monoid $(\F_q,\times)$.
Then $\# X_\Z(\F_q)=\# X(D_q)$, where $X(D)=\Hom(D,X)$ as usual.
For an integer $k\ge 2$ let $C_{k-1}$ denote the cyclic group of $k-1$ elements and let $D_k$ be the monoid $C_{k-1}\cup \{0\}$, where $x\cdot 0=0$. 
Note that if $q$ is a prime power, then $D_q\cong (\F_q,\times)$, where $\F_q$ is the field of $q$ elements.

Fix a covering of $X$ by affines $U_i=\spec A_i$. Since $\spec(D_k)$ consists of two points, the generic, which always maps to the generic point and the closed point, it follows that
$$
X(\spec(D_k))\= \bigcup_{i} U_i(\spec(D_k)),
$$
and thus the cardinality of the right hand side may be written as an alternating sum of terms of the form
$$
\# U_{i_1}\cap\dots\cap U_{i_s}(\spec(D_k)).
$$
Now $U_{i_1}\cap\dots\cap U_{i_s}$ is itself a union of affines and so this term again becomes an alternating sum of similar terms.
This process stops as $X$ is a finite set.
Therefore, to prove the theorem, it suffices to assume that $X$ is affine.

So we assume that $X=\spec(A)$ for a finitely generated monoid $A$.
In this case $X(\spec(D_k))=\Hom(A,D_k)$.
For a given monoid morphism $\ph : A\to D_k$ we have that $\ph^{-1}(\{ 0\})$ is a prime ideal in $A$, call it $\p$.
Then $\ph$ maps $S_\p= A\setminus \p$ to the group $C_{k-1}$.
So $\Hom(A,D_k)$ may be identified with the disjoint union of the sets $\Hom(S_\p,C_{k-1})$ where $\p$ ranges over $\spec(A)$.
Now $C_{k-1}$ is a group, so every homomorphism from $S_\p$ to $C_{k-1}$ factorises over the quotient group $\Quot(S_\p)$ and one gets $\Hom(S_\p,C_{k-1})=\Hom(\Quot(S_\p),C_{k-1})$.
Note that $\Quot(S_\p)$ is the group of units in the stalk $\CO_{X,\p}$ of the structure sheaf, therefore does not depend on the choice of the affine neighbourhood.
The group $\Quot(S_\p)$  is a finitely generated abelian group.
Let $r$ be its rank and $e$ its exponent.
If $e$ is coprime to $k-1$, then there is no non-trivial homomorphism from the torsion part of $\Quot(S_\p)$ to $C_{k-1}$ and so in that case $\#\Hom(S_p, C_{k-1})=(k-1)^r$.
This proves the existence of $e$ and $N$ and finishes the proof of Theorem \ref{zetapoly}.
\qed

{\bf Remark 1.}
We have indeed proven more than Theorem \ref{zetapoly}.
For an $\F_1$-scheme $X$ of finite type we define $X(\F_q)=\Hom(\spec(\F_q),X)$, where the $\Hom$ takes place in the category of $\F_1$-schemes, and $\F_q$ stands for the multiplicative monoid of the finite field.
It follows that
$$
X(\F_q)\ \cong\ X_\Z(\F_q).
$$
Further, for $k\in\N$ one sets $\F_k=D_k$ then this notation is consistent and we have proven above,
$$
(k-1,e)\=1\ \ \Rightarrow\ \ \# X(\F_k)\= N(k),
$$
where $e$ now is a well defined number, the exponent of $X$.
Further it follows from the proof, that the degree of $N$ is at most equal to the rank of $X$, which is defined as the maximum of the ranks of the local monoids $\CO_\p$, for $\p\in X$.

{\bf Remark 2.}
As the proof of Theorem \ref{zetapoly} shows, the zeta-polynomial $N_X$ of $X$, does actually not depend on the structure sheaf $\CO_X$, but on the subsheaf of units $\CO_X^\times$, where for every open set $U$ in $X$ the set $\CO_X^\times(U)$ is defined to be the set of sections $s\in\CO_X(U)$ such that $s(p)$ lies in $\CO_{X,p}^\times$ for every $p\in U$.
We therefore call $\CO_X^\times$ the \emph{zeta sheaf} of $X$.

\section{K-theory}

In this section we give two definitions of K-theory over $\F_1$ and we show that they do coincide for groups, but not in general.
This approach follows Quillen \cite{Qu}.

\subsection{The $+$-construction}
Let $A$ be a monoid.
Recall from \cite{F1} that $\GL_n(A)$ is the group of all $n\times n$ matrices with exactly one non-zero entry in each row and each column, and this entry being an element of the unit group $A^\times$. We also write $A^\times$ as the stalk $A_c$ at the closed point $c$ of $\spec(A)$.
In other words, we have
$$
\GL_n(A)\ \cong\   A_c^n\rtimes \Per(n),
$$
where $\Per(n)$ is the permutation group in $n$ letters, acting on $A_c^n$ by permuting the co-ordinates.

There is a natural embedding $\GL_n(A)\hookrightarrow \GL_{n+1}(A)$ by setting the last co-ordinate equal to $1$.
We define the group
$$
\GL(A)\df \lim_{\stackrel{\longrightarrow}{n}}\GL_n(A).
$$
Similar to the K-theory of rings \cite{Qu} for $j\ge 0$ we define
$$
K_j^+(A)\df \pi_j(\BGL(A)^+),
$$
where $\BGL(A)$ is the classifying space of $\GL(A)$, the $+$ signifies the $+$-construction, and $\pi_j$ is the $j$-th homotopy group.
For instance, $K_j^+(\F_1)$ is the $j$-th stable homotopy group of the spheres \cite{Priddy}.

\subsection{The Q-construction}
A category is called \emph{balanced},if every morphism which is epi and mono, already has an inverse, i.e., is an isomorphism.

Let $\CC$ be a category.
An object $I\in \CC$ is called \emph{injective} if for every monomorphism $M\hookrightarrow N$ the induced map $\Mor(N,I)\to \Mor(M,I)$ is surjective.
Conversely, an object $P\in\CC$ is called \emph{projective} if for every epimorphism $M\onto N$ the induced map $\Mor(P,M)\to\Mor(P,N)$ is surjective.
We say that $\CC$ \emph{has enough injectives} if for every $A\in\CC$ there exists a monomorphism $A\hookrightarrow I$, where $I$ is an injective object.
Likewise, we say that $\CC$ \emph{has enough projectives} if for every $A\in\CC$ there is an epimorphism $P\onto A$ with $P$ projective.

A category $\CC$ is \emph{pointed} if it has an object $0$ such that for every object $X$ the sets $\Mor(X,0)$ and $\Mor(0,X)$ have exactly one element each.
The zero object is uniquely determined up to unique isomorphy.
In every set $\Mor(X,Y)$ there exists a unique morphism which factorises over the zero object, this is called the zero morphism.
In a pointed category it makes sense to speak of kernels and cokernels.
Kernels are always mono and cokernels are always epimorphisms.
A sequence
$$
\begin{diagram}
\node 0\arrow{e}\node{X}\arrow{e,t}{i}\node{Y}\arrow{e,t}{j}\node{Z}\arrow{e}\node 0
\end{diagram}
$$
is called \emph{strong exact}, if $i$ is the kernel of $j$ and $j$ is the cokernel of $i$.
We say that the sequence \emph{splits}, if it is isomorphic to the natural sequence
$$
0\ \to\ X\ \to\ X\oplus Z\ \to\  Z\ \to\ 0.
$$

Assume that kernels and cokernels always exist. Then every kernel is the kernel of its cokernel and every cokernel is the cokernel of its kernel.
For a morphism $f$ let $\im(f)=\ker(\coker(f))$ and $\coim(f)=\coker(\ker(f))$.
If $\CC$ has enough projectives, then the canonical map $\im(f)\to\coim(f)$ has zero kernel and if $\CC$ has enough injectives, then this map has zero cokernel.

Let $\CC$ be a pointed category and $\CE$ a class of strong exact sequences.
The class $\CE$ is called \emph{closed under isomorphism}, or simply \emph{closed} if every seqeunce isomorphic to one in $\CE$, lies in $\CE$.
Every morphism occurring in a sequence in $\CE$ is called an $\CE$-morphism.

A balanced pointed category $\CC$, together with a closed class
$\CE$ of strong exact sequences  is called a \emph{quasi-exact category} if 
\begin{itemize}
\item for any two objects $X,Y$ the natural sequence
$$
0\to X\to X\oplus Y\to Y\to 0
$$
belongs to $\CE$,
\item
the class of $\CE$-kernels is closed under composition and base-change by $\CE$-cokernels,
likewise, the class of $\CE$-cokernels is closed under composition and base change by $\CE$-kernels.
\end{itemize}

Let $(\CC,\CE)$ be a quasi-exact category.
We define the category $Q\CC$ to have the same objects as $\CC$, but a morphism from $X$ to $Y$ in $Q\CC$ is an isomorphism  class of diagrams of the form
$$
\begin{diagram}
\node S\arrow{e,J}\arrow{s,A}\node Y\\
\node {X,}
\end{diagram}
$$
where the horizontal map is a $\CE$-kernel in $\CC$ and the vertical map is a $\CE$-cokernel.
The composition of two Q-morphisms
$$
\begin{diagram}
\node S\arrow{e,J}\arrow{s,A}\node Y\\
\node {X,}
\end{diagram}
\qquad
\begin{diagram}
\node T\arrow{e,J}\arrow{s,A}\node Z\\
\node {Y,}
\end{diagram}
$$
is given by the base change $S\times_YT$ as follows,
$$
\begin{diagram}
\node{S\times_Y T}\arrow{e,J}\arrow{s,A}
	\node{T}\arrow{e,J}\arrow{s,A}
		\node{Z}\\
\node S\arrow{e,J}\arrow{s,A}\node Y\\
\node {X.}
\end{diagram}
$$

Every $\CE$-kernel $\begin{diagram}
\node{X}\arrow{e,J,t}{i:}\node Y
\end{diagram}$ gives rise to a morphism $i_!$ in $Q\CC$, and every $\CE$-cokernel $\begin{diagram}
\node Z\arrow{e,A,t}{p:}\node Z
\end{diagram}$ gives rise to a morphism $p^! : X\to Z$ in $Q\CC$.
By definition, every morphism in $Q\CC$ factorises as $i_!p^!$ uniquely up to isomorphism.

Let $(\CC,\CE)$ be a small quasi-exact category.
Then the classifying space $BQ\CC$ is defined.
Note that for every object $X$ in $Q\CC$ there is a morphism from $0$ to $X$, so that $BQ\CC$ is path-connected.
We consider the fundamental group $\pi_1(BQ\CC)$ as based at a zero $0$ of $\CC$.

\begin{theorem}
The fundamental group $\pi_1(BQ\CC)$ is canonically isomorphic to the Grothendieck group $K_0(\CC)=K_0(\CC,\CE)$.
\end{theorem}

\prf
This proof is taken from \cite{Qu}, where it is done for exact categories, we repeat it for the convenience of the reader.
The Grothendieck group $K_0(\CC,\CE)$ is the abelian group with one generator $[X]$ for each object $X$ of $\CC$ and a relation $[X=[Y][Z]$ for every strong exact sequence
$$
\begin{diagram}
\node 0\arrow{e}\node{Y}\arrow{e,J}\node{X}\arrow{e,A}\node{Z}\arrow{e}\node{0}
\end{diagram}
$$
in $\CE$.
According to Proposition 1 of \cite{Qu}, it suffices to show that for a morphism-inverting functor 
$F: Q\CC\to \Sets$ the group $K_0(\CC)$ acts 
naturally on $F(0)$ and that the resulting 
functor from the category $\CF$ of all such $F$ to $K_0(\CC)$-sets is an equivalence of categories.

For $X\in\CC$ let $i_X$ denote the zero kernel $0\to X$,  and let $j_X$ be the zero cokernel $X\to 0$.
Let $\CF'$ be the full subcategory of $\CF$ consisting of all $F$ such that $F(X)=F(0)$ and $F(i_{X!})=\id_{F(0)}$ for every $X$.
Any $F\in\CF$ is isomorphic to an object of $\CF'$, so it suffices to show that $\CF'$ is equivalent to $K_0(\CC)$-sets.
So let $F\in\CF'$, for a kernel $I:X\hookrightarrow Y$ we have $ii_X=i_Y$, so that $F(i_!)=\id_{F(0)}$.
Given a strong exact sequence
$$
\begin{diagram}
\node{0}\arrow{e}\node{X}\arrow{e,t,J}{i}\node{Y}
\arrow{e,t,A}{j}\node{Z}\arrow{e}\node{0,}\\
\end{diagram}
$$
we have $j^!i_{Z!}=i_!j_X^!$, hence $F(j^!)=F(j_X^!)\in\Aut(F(0))$.
Also,
$$
F(j_Y^!)\= F(j^!j_Z^!)\= F(j_X^!)F(j_Z^!).
$$
So by the universal property of $K_0(\CC)$, there is a unique homomorphism from $K_0(\CC)$ to $\Aut(F(0))$ such that $[X]\mapsto F(j_X^!)$.
So we have a natural action of $K_0(\CC)$ on $F(0)$, hence a functor from $\CF'$ to $K_0(\CC)$-sets given by $F\mapsto F(0)$.

The other way round let $S$ be a $K_0(\CC)$-set, and let $F_S: \QB\to\Sets$ be the functor defined by $F_S(X)=S$, $F_S(i_!j^!)=$ multiplication by $[\ker j]$ on $S$.
To see that this is indeed a functor, it suffices to show that $F_S(j^!i_!)=F_S(j^!)$.
It holds $j^!i_!=i_{1!}j_1^!$, where $i_1$ and $j_1$ are given by the cartesian diagram
$$
\begin{diagram}
\node{A}\arrow{e,t,J}{i_1}\arrow{s,l,A}{j_1}
	\node{X}\arrow{s,r,A}{j}\\
\node{Z}\arrow{e,t,J}{i}
	\node{Y.}
\end{diagram}
$$
It follows $F_S(j^!i_!)=F_S(i_{1!}j_1^!)=[\ker j_1]$.
Using the cartesian diagram one sees that $\ker j_1$ is isomorphic to $\ker j$.
It is easy to verify that the two functors given are inverse to each other up to isomorphism, whence the theorem.
\qed

This theorem motivates the following definition,
$$
K_i(\CC,\CE)\df\pi_{i+1}(BQ\CC).
$$
For a monoid $A$ we let $\CP$ be the category of finitely generated pointed projective $A$-modules, or rather a small category equivalent to it, and we set
$$
K_i^Q(A)\df K_i(\CP,\CE),
$$ where $\CE$ is the class of sequences in $\CP$ which are strong exact in the category of all modules.
These sequences all split, which establishes the axioms for a quasi-exact category.

The two $K$-theories we have defined, do not coincide.
For instance for the monoid of one generator $A=\{1,a\}$ with $a^2=a$ one has
$$
K_0^+(A)\= \Z,\qquad K_0^Q(A)\= \Z\times\Z.
$$
The reason for this discrepancy is that $K_i^+(A)$ only depends on the group of units $A^\times$, but $K_i^Q(A)$ is sensible to the whole structure of $A$.
So these two K-theories are unlikely to coincide except when $A$ is a group, in which case they do, as the last theorem of this paper shows,

\begin{theorem}
If $A$ is an abelian group, then $K_i^+(A)=K_i^Q(A)$ for every $i\ge 0$.
\end{theorem}

\prf For a group each projective module is free, hence the proof of Grayson \cite{Grays} of the corresponding fact for rings goes through.
\qed

So, if $A$ is a group, this defines $K_i(A)$ unambiguously.
In particular, computations of Priddy \cite{Priddy} show that $K_i(\F_1)=\pi_s^i$ is the $i$-th stable homotopy group of the spheres. 
Based on this, one can use the Q-construction to show that if $A$ is an abelian group, then
$$
K_i(A)\= \begin{cases} \Z\times A & i=0,\\ \pi_i^s & i>0.
\end{cases}
$$
For an arbitrary monoid $A$ we conclude that $K_i^+(A)=K_i^+(A^\times)=K_i(A^\times)$, which we now can express in terms of the stable homotopy groups $\pi_i^s$.

Further, for every $A$ one has a canonical homomorphism $K_i^+(A)\to K_i^Q(A)$ given by the map $K^Q(A^\times)\to K^Q(A)$.
The latter comes about by the fact that every projective $A^\times$-module is free.
Note that general functoriality under monoid homomorphism is granted for $K^+$, but not for $K^Q$.
This contrasts the situation of rings, and has its reason in the fact that not every projective is a direct summand of a free module.

{\small Mathematisches Institut\\
Auf der Morgenstelle 10\\
72076 T\"ubingen\\
Germany\\
\tt deitmar@uni-tuebingen.de}
\end{multicols}
\end{document}